\documentclass[a4paper]{amsart}
\usepackage[english]{babel}

\usepackage{amsmath,amssymb,amsthm}
\theoremstyle{plain}
\newtheorem{theorem}{Theorem}[section]
\newtheorem{corollary}[theorem]{Corollary}
\newtheorem{prop}[theorem]{Proposition}
\newtheorem{lemma}[theorem]{Lemma}

\theoremstyle{definition}
\newtheorem{remark}[theorem]{Remark}

\newcommand{\ec}{\ensuremath{\mathrm{co}}}
\newcommand{\ext}[1][X^*]{\ensuremath{\mathrm{ext}\left(B_{#1}\right)}}
\newcommand{\extr}{\ensuremath{\mathrm{ext}}}
\DeclareMathOperator{\re}{Re\,}

\newcounter{equi1}
\newenvironment{equi}
{\begin{list}
        {(\roman{equi1})}
        {\setlength{\itemsep}{0ex plus 0.2ex minus 0ex}
         \setlength{\topsep}{0ex}
         \setlength{\parsep}{0ex}
         \setlength{\labelwidth}{7ex}
         \usecounter{equi1}}}
{\end{list}}

\begin{document}

\title[The Alternative Daugavet Property]{The Alternative Daugavet Property
of $\mathbf{C^*}$-algebras and $\mathbf{JB^*}$-triples}
\thanks{Research partially supported by Spanish MCYT project no.\ BFM2003-01681
and Junta de Andaluc\'{\i}a grant FQM-185}
 \date{November 22th, 2004}

\maketitle

\centerline{\textsc{\large Miguel Mart\'{\i}n}}

\begin{center} Departamento de An\'{a}lisis Matem\'{a}tico \\ Facultad de
Ciencias \\ Universidad de Granada \\ 18071 Granada, SPAIN \\
\emph{Email address:} \texttt{mmartins@ugr.es} \end{center}

 \thispagestyle{empty}

\begin{abstract}
A Banach space $X$ is said to have the alternative Daugavet
property if for every (bounded and linear) rank-one operator
$T:X\longrightarrow X$ there exists a modulus one scalar $\omega$
such that $\|Id + \omega T\|= 1 + \|T\|$. We give geometric
characterizations of this property in the setting of
$C^*$-algebras, $JB^*$-triples and their isometric preduals.
\end{abstract}

\maketitle

\vspace{1cm}

For a real or complex Banach space $X$, we write $X^*$ for its
topological dual and $L(X)$ for the Banach algebra of bounded
linear operators on $X$. We denote by $\mathbb{T}$ the set of
modulus-one scalars.

Following \cite{MaOi}, we say that a Banach space $X$ has the
\emph{alternative Daugavet property} if every rank-one operator
$T\in L(X)$ satisfies the norm identity
\begin{equation*}\label{aDE}\tag{\textrm{aDE}}
\max_{\omega\in\mathbb{T}}\|Id +\omega\, T\| = 1 + \|T\|.
\end{equation*}
In such a case, all weakly compact operators on $X$ also satisfy
Eq.~\eqref{aDE} (see \cite[Theorem~2.2]{MaOi}). It is clear that a
Banach space $X$ has the alternative Daugavet property whenever
$X^*$ has, but it is shown in \cite[Example~4.4]{MaOi} that the
reverse result does not hold.

Observe that Eq.~\eqref{aDE} for an operator $T$ just means that
there exists a modulus-one scalar $\omega$ such that the operator
$S=\omega T$ satisfies the (usual) \emph{Daugavet equation}
\begin{equation*}\label{DE}\tag{\textrm{DE}}
\|Id + S\| = 1 + \|S\|.
\end{equation*}
Therefore, the \emph{Daugavet property} (i.e., every rank-one
--equivalently, every weakly compact-- operator satisfies
Eq.~\eqref{DE} \cite[Theorem~2.3]{KSSW}) implies the alternative
Daugavet property. Examples of spaces having the Daugavet property
are $C(K)$ and $L_1(\mu)$, provided that $K$ is perfect and $\mu$
does not have any atoms (see \cite{Wer0} for an elementary
approach), and certain function algebras such as the disk algebra
$A(\mathbb{D})$ or the algebra of bounded analytic functions
$H^\infty$ \cite{WerJFA,Woj}. The state-of-the-art on the subject
can be found in \cite{KSSW,WerSur}. For very recent results we refer
the reader to \cite{BM,BKSW,KaKaWe,KaWe,Oik} and references therein.

Although the definition of the alternative Daugavet property is
very recent, Eq.~\eqref{aDE} appears explicitly in several papers
from the 80's as \cite{Abra0,Abra,Holub,Holub2,Wer0}, where it is
proved that this equation is satisfied by all operators $T\in
L(X)$ whenever $X=C(K)$ or $X=L_1(\mu)$. Actually, this result
appeared in the 1970 paper \cite{D-Mc-P-W}, where Eq.~\eqref{aDE}
is related to a constant introduced by G.~Lumer in 1968, the
numerical index of a Banach space. Let us give the necessary
definitions. Given an operator $T\in L(X)$, the \emph{numerical
radius} of $T$ is
$$v(T):=\sup\{|x^*(T x)|\ : \ x\in X,\ x^*\in X^*,\
\|x^*\|=\|x\|=x^*(x)=1\}.$$ The \emph{numerical index} of the
space $X$ is $$n(X):=\max \left\{k\geqslant 0 \ : \
k\|T\|\leqslant v(T)\ \text{for all }T\in L(X)\right\}.$$ It was
shown in \cite[pp.~483]{D-Mc-P-W} that
\begin{center}
an operator $T\in L(X)$ satisfies Eq.~\eqref{aDE} if and only if
$v(T)=\|T\|$,
\end{center}
thus, $n(X)=1$ if and only if Eq.~\eqref{aDE} is satisfied by all
operators in $L(X)$. Therefore, the numerical index~$1$ implies
the alternative Daugavet property. Examples of spaces with
numerical index~$1$ are $L_1(\mu)$ and its isometric preduals, for
any positive measure $\mu$. For more information and background,
we refer the reader to the monographs by F.~Bonsall and J.~Duncan
\cite{B-D1,B-D2} and to the survey paper \cite{Mar}. Recent
results can be found in \cite{KMR,LMP,MaMeRo,M-V,Oik04} and
references therein. Let just mention that one has $v(T^*)=v(T)$
for every $T\in L(X)$, where $T^*$ is the adjoint operator of $T$
(see \cite[\S 9]{B-D1}) and it clearly follows that
$n(X^*)\leqslant n(X)$ for every Banach space $X$. The question if
this is actually an equality seems to be open.

The alternative Daugavet property does not imply neither the
Daugavet property nor numerical index $1$. Indeed, {\slshape the
space $c_0\oplus_\infty C([0,1],\ell_2)$ has the alternative
Daugavet property, but it does not have the Daugavet property and
it does not have numerical index $1$} (see
\cite[Example~3.2]{MaOi} for the details).

In \cite{MaOi}, the authors characterize the $JB^*$-triples having
the alternative Daugavet property as those whose minimal
tripotents are diagonalizing. It is also proved there that the
predual of a $JBW^*$-triple has the alternative Daugavet property
if and only if the triple does. The necessary definitions and
basic results on $JB^*$-triples are presented in
section~\ref{sec-JB}.

In the present paper we give geometric characterizations of the
alternative Daugavet property for $JB^*$-triples and for their
isometric preduals. In particular, our results contain the above
mentioned algebraic characterizations given in \cite{MaOi}, but
our proofs are independent. To state the main results of the paper
we need to fix notation and recall some definitions.

Let $X$ be a Banach space. The symbols $B_X$ and $S_X$ denote,
respectively, the closed unit ball and the unit sphere of $X$. We
write $\ec(A)$ for the convex hull of the set $A$, and we will
denote by $\extr(B)$ the set of extreme points in the convex set
$B$. Let us fix $u$ in $S_X$. We define the set $D(X,u)$ of all
\emph{states} of $X$ relative to $u$ by $$D(X,u):= \{f\in B_{X^*}\
: \ f(u)=1\},$$ which is a non-empty $w^*$-closed face of
$B_{X^*}$. The norm of $X$ is said to be \emph{Fr\'{e}chet-smooth} or
\emph{Fr\'{e}chet differentiable} at $u\in S_X$ whenever there exists
$\displaystyle \lim _{\alpha \rightarrow 0} \frac{\Vert u+\alpha
x\Vert -1}{\alpha}$ uniformly for $x\in B_X$. A point $x\in S_X$
is said to be an \emph{strongly exposed point} if there exists
$f\in D(X,x)$ such that $\lim \|x_n - x\|=0$ for every sequence
$(x_n)$ of elements of $B_{X}$ such that $\lim\,\re f(x_n)=1$. If
$X$ is a dual space and $f$ is taken from the predual, we say that
$x$ is a \emph{$w^*$-strongly exposed point}. It is known that $x$
is strongly exposed if and only if there is a point of
Fr\'{e}chet-smoothness in $D(X,x)$, which is actually an strongly
exposing functional for $x$  (see \cite[Corollary~I.1.5]{DGZ}).
Finally, if $X$ and $Y$ are Banach spaces, we write $X\oplus_1 Y$
and $X\oplus_\infty Y$ to denote, respectively, the $\ell_1$-sum
and the $\ell_\infty$-sum of $X$ and $Y$.

The main results of the paper are the characterizations of the
alternative Daugavet property for $JB^*$-triples and preduals of
$JBW^*$-triples given in Theorems \ref{th:comm+DPforC} and
\ref{theorem-JBW} respectively. For a $JB^*$-triple $X$, the
following are equivalent:
\begin{enumerate}
\item[$(i)$] $X$ has the alternative Daugavet property.
\item[$(ii)$] $|x^{**}(x^*)|=1$ for every $x^{**}\in \ext[X^{**}]$ and every
$w^*$-strongly exposed point $x^*$ of $B_{X^*}$.
\item[$(iii)$] Each elementary triple ideal of $X$ is equal
to $\mathbb{C}$.
\item[$(iv)$] There exists a closed triple ideal $Y$ with
$c_0(\Gamma)\subseteq Y \subseteq \ell_\infty(\Gamma)$ for
convenient index set $\Gamma$, and such that $X/Y$ is non-atomic.
\item[$(v)$] All minimal tripotents of $X$ are diagonalizing.
\end{enumerate}

For the predual $X_*$ of a $JBW^*$-triple $X$, the following are
equivalent:
\begin{enumerate}
\item[$(i)$] $X_*$ has the alternative Daugavet property,
\item[$(ii)$] $X$ has the alternative Daugavet property,
\item[$(iii)$] $|x(x_*)|=1$ for every $x\in \ext[X]$ and every $x_*\in
\ext[X_*]$.
\item[$(iv)$] $X=\ell_\infty(\Gamma)\oplus_\infty \mathcal{N}$ for
suitable index set $\Gamma$ and non-atomic ideal $\mathcal{N}$.
\end{enumerate}

Let us mention that geometric characterizations of the Daugavet
property for $C^*$-algebras, $JB^*$-triples, and their isometric
preduals, can be found in \cite{BM}.

The outline of the paper is as follows. In section~2 we include
some preliminary results on the alternative Daugavet property and
on Banach spaces with numerical index~$1$. Section~3 is devoted to
the above cited characterizations of the alternative Daugavet
property for $JB^*$-triples and their isometric preduals, and we
dedicate section~4 to particularize these result to $C^*$-algebras
and von Neumann preduals, and to compile some results on
$C^*$-algebras and von Neumann preduals with numerical index $1$.

\section{Preliminary results}
We devote this section to summarize some results concerning the
alternative Daugavet property and Banach spaces with numerical
index~$1$, which we will use along the paper. Most of them appear
implicity in \cite{LMP,MaOi}, but there are no explicit
references. Therefore, we include them here for the sake of
completeness.

Our starting point is McGregor's characterization of
finite-dimensional spaces with numerical index $1$ \cite[Theorem
3.1]{Mc}: \emph{A finite-dimensional space $X$ satisfies $n(X)=1$
if and only if $|x^*(x)|=1$ for all extreme points $x\in B_X$ and
$x^*\in B_{X^*}$}. For infinite-dimensional $X$, $\ext[X]$ may be
empty (e.g. $c_0$), so the right statement of McGregor's condition
in this case should read
\begin{equation}\label{eq:extr-extr}
|x^{**}(x^*)|=1 \qquad \text{for every $x^*\in \ext[X^*]$ and
every $x^{**}\in \ext[X^{**}]$.}
\end{equation}
One can easily show that this condition is sufficient to ensure
$n(X)=1$. Actually, a bit more general result can be also proved
easily.

\begin{lemma}\label{lemma-sufficient-for-nx=1}
Let $X$ be a Banach space and let $A$ be a subset of $S_{X^*}$
such that $B_{X^*}=\overline{\ec(A)}^{w^*}$. If
$$
|x^{**}(a^*)|=1 \qquad \text{for every $a^*\in A$ and every
$x^{**}\in \ext[X^{**}]$},
$$
then $n(X)=1$. In particular,
condition~\eqref{eq:extr-extr} implies $n(X)=1$.
\end{lemma}

\begin{proof}
Fix $T\in L(X)$. Since $T^*$ is $w^*$-continuous, for every
$\varepsilon>0$ we can find $a^*\in A$ such that
$\|T^*(a^*)\|\geqslant \|T\|-\varepsilon$. Now, take $x^{**}\in
\ext[X^{**}]$ with
$$|x^{**}(T^*(a^*))|=\|T^*(a^*)\|,$$ and observe that $|x^{**}(a^*)|=1$,
giving
\begin{equation*}
v(T)=v(T^*)\geqslant |x^{**}(T^*(a^*))|\geqslant
\|T\|-\varepsilon. \qedhere
\end{equation*}
\end{proof}

We do not know if condition \eqref{eq:extr-extr} is also necessary
in order to have numerical index~$1$, but considering strongly
exposed points, in \cite[Lemma~1]{LMP} a (weaker) necessary
condition is established. Actually, as is noted in the Remark~6 of
the same paper, we can replace the hypothesis of having numerical
index $1$ by the hypothesis of having the alternative Daugavet
property.

\begin{lemma}\label{lemma-LMP}
Let $X$ be a Banach space with the alternative Daugavet property.
Then:
\begin{enumerate}
\item[(a)] $|x^{**}(x^*)|=1$ for every $x^{**}\in \ext[X^{**}]$ and every
$w^*$-strongly exposed point $x^*\in B_{X^*}$.
\item[(b)] $|x^*(x)|=1$ for every $x^*\in \ext[X^*]$ and every strongly exposed
point $x\in B_X$.
\end{enumerate}
\end{lemma}

Recall that the dual unit ball of an Asplund space is the
$w^*$-closed convex hull of its $w^*$-strongly exposed points (see
\cite{Ph}, for instance) so, by just applying the above two
lemmas, we easily get the following characterization.

\begin{corollary}\label{coro-char-Asplund}
Let $X$ be an Asplund space. Then, the following are equivalent:
\begin{equi}
\item $n(X)=1$.
\item $X$ has the alternative Daugavet property.
\item $|x^{**}(x^*)|=1$ for every $x^{**}\in \ext[X^{**}]$ and every
$w^*$-strongly exposed point $x^*\in B_{X^*}$.
\end{equi}
\end{corollary}

For similar results concerning Banach spaces with the
Radon-Nikod\'{y}m property and numerical index $1$ we refer the reader
to \cite{MarRNP}.

Finally, let us mention that we do not know of any general
characterization of Banach spaces having numerical index $1$ which
does not involve operators. In particular, we do not know if
condition \eqref{eq:extr-extr} is such a characterization.

\section{$JB^*$-triples and preduals of $JBW^*$-triples}\label{sec-JB}

We recall that a \emph{$JB^*$-triple} is a complex Banach space
$X$ with a continuous triple product $\{\cdots\}:\, X\times
X\times X\longrightarrow X$ which is linear and symmetric in the
outer variables, conjugate-linear in the middle variable, and
satisfies:
\begin{enumerate}
\item For all $x$ in $X$, the mapping $y\longmapsto \{xxy\}$ from
$X$ to $X$ is a hermitian operator on $X$ and has nonnegative
spectrum.
\item The \emph{main identity} $$\{ab\{xyz\}\}=\{\{abx\}yz\}-
\{x\{bay\}z\}+\{xy\{abz\}\} $$ holds for all $a,b,x,y,z$ in $X$.
\item $\Vert \{xxx\}\Vert =\Vert x\Vert ^{3}$ for every $x$ in
$X$.
\end{enumerate}
Concerning Condition (1) above, we also recall  that  a  bounded
linear operator $T$ on a complex Banach space $X$ is said to be
\emph{hermitian} if $\Vert \exp (irT)\Vert =1$ for every
$r\in\mathbb{R}$. The main interest of $JB^*$-triples relies on
the fact that, up to biholomorphic equivalence, there are no
bounded symmetric domains in complex Banach spaces others than the
open unit balls of $JB^*$-triples (see \cite{K,K1}).

Every $C^*$-algebra becomes a $JB^*$-triple under the triple
product $$\{xyz\}:=\frac{1}{2}(xy^*z+zy^*x).$$ Moreover,
norm-closed subspaces of $C^*$-algebras closed under the above
triple product form the class of $JB^*$-triples known as
\emph{$JC^*$-triples}.

By a \emph{$JBW^*$-triple} we mean a $JB^*$-triple whose
underlying Banach space is a dual space in metric sense. It is
known (see \cite{BaTi}) that every $JBW^*$-triple has a unique
predual up to isometric linear isomorphisms and its triple product
is separately $w^*$-continuous in each variable. We will apply
without notice that the bidual of every $JB^*$-triple $X$ is a
$JBW^*$-triple under a suitable triple product which extends the
one of $X$ \cite{Dineen}.

Let $X$ be a $JB^*$-triple. Given $a,b\in X$, we write $D(a,b)$
for the multiplication operator $x\longmapsto \{abx\}$. The
elements $a$ and $b$ are said to be \emph{orthogonal} if
$D(a,b)=0$ (equivalently, $D(b,a)=0$). For $a\in X$, the conjugate
linear operator $x\longmapsto \{axa\}$ is denoted by $Q_a$. An
element $u\in X$ is said to be a \emph{tripotent} if $\{uuu\}=u$.
Associated to any tripotent $u$, we define the \emph{Peirce
projections}
$$
P_2(u) = Q_u^2, \qquad P_1(u) = 2(D(u,u) - Q_u^2), \qquad P_0(u) =
Id - 2D(u,u) + Q_u^2,
$$
which are mutually orthogonal with sum $Id$, and ranges
$$X_j(u):=P_j(u)(X)=\Big\{x\ : \ \{uux\}=\frac{j}{2}x\Big\}\qquad
(j = 0,1,2),$$ giving $X=X_2(u)\oplus X_1(u)\oplus X_0(u)$
(\emph{Peirce decomposition}). A tripotent $u\in X$ is said to be
\emph{minimal} if $u\neq 0$ and $X_2(u)=\mathbb{C} u$
(equivalently, $\{uXu\}=\mathbb{C} u$) and $u$ is said to be
\emph{diagonalizing} if $X_1(u)=0$. A \emph{triple ideal} of $X$
is a subspace $M$ of $X$ such that $\{XXM\}+ \{XMX\}\subseteq M$.
It is well-known that the triple ideals of a $JB^*$-triple are
precisely its $M$-ideals \cite[Theorem~3.2]{BaTi}. The
\emph{socle} of $X$, $K_0(X)$, is the closed linear span of the
minimal tripotents of $X$. Then, $K_0(X)$ is a triple ideal of $X$
(the one generated by the minimal tripotents in $X$), which is
equals to the $c_0$-sum of every elementary triple ideal of $X$
\cite[Lemma~3.3]{BC92Pacific}. $X$ is said to be \emph{non-atomic}
if $K_0(X)=0$, that is, if it does not contain any minimal
tripotent. By \cite[Theorem~4.7]{MaOi} or
\cite[Corollary~3.4]{BM}, a $JB^*$-triple is non-atomic if and
only if it has the Daugavet property. On the other side, a
$JBW^*$-tripe $Y$ is said to be \emph{atomic} if it equals the
$w^*$-closure of it socle. The $JBW^*$-triple $Y$ is said to be a
\emph{factor} if it cannot be written as an $\ell_\infty$-sum of
two nonzero ideals, or equivalently \cite{Horn-predual} if $0$ and
$Y$ are the only $w^*$-closed ideals in $Y$. An specially
important class of $JBW^*$-triples are the so-called \emph{Cartan
factors}. This class falls into six subclasses as follows, where
$H$ and $H'$ are arbitrary complex Hilbert spaces and
$J:H\longrightarrow H$ is a conjugation (that is, a conjugate
linear isometry with $J^2=Id$):
\begin{enumerate}
\item[(1)] \emph{rectangular:} $L(H,H')$;
\item[(2)] \emph{symplectic:} $\{x\in L(H)\ : \ Jx^*J=-x\}$;
\item[(3)] \emph{hermitian:} $\{x\in L(H)\ : \ Jx^*J=x\}$;
\item[(4)] \emph{spin:} $H$ with $\dim(H)\geqslant 3$, with convenient
(non-Hilbertian) equivalent norm.
\item[(5)] $M_{1,2}(\mathbb{O})$: the $1\times 2$ matrices over the complex
octonians $\mathbb{O}$;
\item[(6)] $H_3(\mathbb{O})$: the hermitian $3\times 3$ matrices over
$\mathbb{O}$.
\end{enumerate}
The rectangular, symplectic and hermitian factors have the usual
$JC^*$-triple product. For the definition of the other triple
products see \cite{K1}. The Cartan factors are appropriately
``irreducible'' $JBW^*$-triples \cite{Horn-classification} in
terms of which there is a useful representation theory. For
instance, every $JB^*$-triple is isometrically isomorphic to a
subtriple of an $\ell_\infty$-sum of Cartan factors \cite{FR86}.

A $JB^*$-triple $M$ is said to be \emph{elementary} if it is
isometric to the socle $K_0(C)$ of a Cartan factor $C$. We have
$K_0(C)^{**}=C$, and that a $JB^*$-triple is elementary if and
only if its bidual is a Cartan factor \cite[Lemma~3.2]{BC91MZ}. Of
course, no elementary $JB^*$-tripe has the Daugavet property. The
following result said that, actually, not so many elementary
$JB^*$-triples have the alternative Daugavet property.

\begin{prop}\label{Cartan}
The only Cartan factor which has the alternative Daugavet property
is $\mathbb{C}$. Actually, the only elementary $JB^*$-triple which
has the alternative Daugavet property is $\mathbb{C}$.
\end{prop}

To prove the above proposition, we need the following easy
consequence of \cite[Proposition~4]{FR85} and
\cite[Corollary~2.11]{BeRo} (see the proof of
\cite[Theorem~3.2]{BM}). For the sake of completeness, we include
a proof here.

\begin{lemma}\label{lemma-extreme-stexposed}
Let $X$ be a $JBW^*$-triple. Then, every extreme point $x_*$ of
$B_{X_*}$ is strongly exposed by a minimal tripotent of $X$.
\end{lemma}

\begin{proof}
Given $x_*\in\ext[X_*]$, \cite[Proposition~4]{FR85} assures the
existence of a minimal tripotent $u$ of $X$ such that $u(x_*)=1$.
But $u$ is a point of Fr\'{e}chet-smoothness of the norm of $X$
\cite[Corollary~2.11]{BeRo} so, as we commented in the
introduction, this implies that $x_*$ is strongly exposed by $u$
(see \cite[Corollary~I.1.5]{DGZ}, for instance).
\end{proof}

\begin{proof}[Proof of Proposition~\ref{Cartan}]
Let $M$ be an elementary $JB^*$-triple and write $C=M^{**}$, which
is a Cartan factor. Suppose that $M$ has the alternative Daugavet
property. By the above lemma, every extreme point $x^*$ of
$B_{M^*}$ is strongly exposed by a minimal tripotent of $M^{**}$,
which is actually in $M$, so $x^*$ is a $w^*$-strongly exposed
point of $B_{M^*}$. Then, Lemma~\ref{lemma-LMP} gives us that
\begin{equation}\label{eq:Cartan}
|e(x_*)|=1 \qquad \ \big(e\in \ext[C],\ x_*\in \ext[C_*]\big).
\end{equation}
We claim that
\begin{equation}\label{eq:Cartan2}
|u(x_*)|\in\{0,1\}
\end{equation}
for every minimal tripotent $u$ of $C$ and every
$x_*\in\ext[C_*]$. Indeed, we fix a minimal tripotent $u$ of $C$
and $x_*\in \ext[C_*]$, and we find a maximal orthogonal family of
minimal tripotents $\{u_i\}_{i\in I}$ such that $u_{i_0}=u$ for
some $i_0\in I$. Now, for each family $\{\lambda_i\}_{i\in I}$ of
elements of $\mathbb{T}$, we consider $e=\sum_{i\in I} \lambda_i
u_i$, which is a complete tripotent of $C$ \cite[Lemma~2.1]{She}
and thus, an extreme point of $B_C$ \cite{KaUp}, so
Eq.~\eqref{eq:Cartan} gives $$\left|\sum_{i\in I} \lambda_i
u_i(x_*)\right| =|e(x_*)|=1.$$ In particular, choosing a suitable
family $\{\lambda_i\}_{i\in I}$, we get
$$\sum_{i\in I}|u_i(x_*)|=1.$$ If we write $\alpha=u_{i_0}(x_*)$,
$\beta =\sum_{i\in I\setminus\{i_0\}} u_i(x_*)$, the above two
equations give that
$$|\alpha|\leqslant 1,\ |\beta|\leqslant 1, \qquad \text{and}\qquad |\alpha + \lambda\beta|=1\ \
\text{for every } \lambda\in \mathbb{T}.$$ This implies that
$|\alpha|,|\beta| \in\{0,1\}$, and the claim is proved.

Now, the set of minimal tripotents of a Cartan factor is
connected; indeed, for Cartan factors of finite-rank this can be
easily deduced from \cite[(4.5) and (4.6)]{Kau}, and the remaining
examples follows by a direct inspection (compare
\cite[Section~3]{Kau}). Therefore, Eq.~\eqref{eq:Cartan2} gives us
that $|u(x_*)|=1$ for every minimal tripotent $u$ of $C$ and every
$x_*\in\ext[C_*]$. This implies that every minimal tripotent of
$C$ is complete (there is no pair of orthogonal minimal tripotents
in $C$) and thus, $\text{rank}(C)=1$ and $C$ is a Hilbert space
\cite{DaFri}; but, obviously, if a Hilbert space satisfies
\eqref{eq:Cartan}, then it is one-dimensional.
\end{proof}

If $X$ is a $JBW^*$-triple, it is well known that $X_*=A\oplus_1
N$, where $A$ is the closed linear span of the extreme points of
$B_{X_*}$, and the unit ball of $N$ has no extreme points
\cite{FR85}. Therefore, $X=\mathcal{A}\oplus_\infty \mathcal{N}$,
where $\mathcal{A}=N^\perp\equiv A^*$ is atomic and
$\mathcal{N}=A^\perp\equiv N^*$ is non-atomic. On one hand, by
\cite[Proposition~ 2.2]{FR86} and
\cite[Corollary~1.8]{Horn-classification}, an atomic
$JBW^*$-triple is the $\ell_\infty$-sum of Cartan factors. On the
other hand, the non-atomic part $\mathcal{N}$ of $X$ has the
Daugavet property, hence the alternative Daugavet property. Now,
since an $\ell_\infty$-sum of Banach spaces has the alternative
Daugavet property if and only if all the summands do
\cite[Proposition~3.1]{MaOi}, the above proposition said that $X$
has the alternative Daugavet property if and only if all the
Cartan factors appearing in the atomic part are equal to
$\mathbb{C}$. This is part of the following result, which
characterizes those $JBW^*$-triples whose preduals have the
alternative Daugavet property.

\begin{theorem}\label{theorem-JBW}
Let $X$ be a $JBW^*$-triple. Then, the following are equivalent:
\begin{equi}
\item $X$ has the alternative Daugavet property.
\item $X_*$ has the alternative Daugavet property.
\item $|x(x_*)|=1$ for every $x\in \ext[X]$ and every $x_*\in
\ext[X_*]$.
\item The atomic part of $X$ is isometrically isomorphic to
$\ell_\infty(\Gamma)$ for convenient set $\Gamma$; i.e.,
$X=\ell_\infty(\Gamma)\oplus_\infty \mathcal{N}$ where
$\mathcal{N}$ is non-atomic.
\end{equi}
\end{theorem}

\begin{proof}
$(i)\Rightarrow(ii)$ is clear.

$(ii)\Rightarrow(iii)$. Just apply Lemmas \ref{lemma-LMP} and
\ref{lemma-extreme-stexposed}.

$(iii)\Rightarrow(iv)$. By the comments preceding this theorem, we
only have to prove that every Cartan factor appearing in the
atomic decomposition of $X$ is equal to $\mathbb{C}$. Let us fix
such a Cartan factor $C$, and let $K$ be an elementary
$JB^*$-triple with $K^{**}=C$. Since $K^*$ is an $L$-summand of
$X_*$ and $K^{**}$ is an $M$-summand of $X$, it is straightforward
to show that condition $(iii)$ implies
\begin{equation*}%\label{eq:thJBW1}
|x^{**}(x^*)|=1 \qquad \big(x^{**}\in \ext[K^{**}],\ x^*\in
\ext[K^*]\big).
\end{equation*}
Now, Lemma~\ref{lemma-sufficient-for-nx=1} implies that $n(K)=1$,
so $K$ has the alternative Daugavet property and
Proposition~\ref{Cartan} gives us that $K$ (and so does $C$) is
equal to $\mathbb{C}$.

$(iv)\Rightarrow(i)$. Since $n\big(\ell_\infty(\Gamma)\big)=1$ and
$\mathcal{N}$ has the Daugavet property, both have the alternative
Daugavet property and so does its $\ell_\infty$-sum
\cite[Proposition~3.1]{MaOi}.
\end{proof}

The equivalences $(i)\Leftrightarrow (ii) \Leftrightarrow (iv)$ of
the above theorem appear in \cite[Theorem~4.6 and
Corollary~4.8]{MaOi} with a different proof.

Just applying Theorem~\ref{theorem-JBW} to the $JBW^*$-triple
$X^{**}$, we get the following result.

\begin{corollary}\label{coro-dualADPimplicanx1}
Let $X$ be a $JB^*$-triple. If $X^*$ has the alternative Daugavet
property, then $|x^{**}(x^*)|=1$ for every $x^{**}\in\ext[X^{**}]$
and every $x^*\in \ext[X^*]$. Thus, $n(X)=1$.
\end{corollary}

\begin{remark}
It is worth mentioning that, for an arbitrary Banach space $Z$, the
condition
$$
|z^*(z)|=1 \qquad \text{for every $z\in\ext[Z]$ and every
$z^*\in\ext[Z^*]$}
$$
does not necessarily imply that $Z$ has the alternative Daugavet
property. Indeed, let $H$ be the 2-dimensional Hilbert space and
let us consider $Z=c_0(H)$. Since $\ext[Z]=\emptyset$, it is clear
that $Z$ satisfies the above condition. But since $Z$ is an
Asplund space and $n(Z)=n(H)<1$ by \cite[Proposition~1]{MarPay},
Corollary~\ref{coro-char-Asplund} gives us that $Z$ does not have
the alternative Daugavet property.
\end{remark}

The last result of the section is a characterization of the
alternative Daugavet property for $JB^*$-triples.

\begin{theorem}\label{th:comm+DPforC}
Let $X$ be a $JB^*$-triple. Then, the following are equivalent:
\begin{equi}
\item $X$ has the alternative Daugavet property.
\item $|x^{**}(x^*)|=1$ for every $x^{**}\in \ext[X^{**}]$ and every
$w^*$-strongly exposed point $x^*$ of $B_{X^*}$.
\item Each elementary triple ideal of $X$ is equal
to $\mathbb{C}$; equivalently, $K_0(X)$ is isometric to
$c_0(\Gamma)$ for some index set $\Gamma$.
\item There exists a closed triple ideal $Y$ with $c_0(\Gamma)\subseteq Y \subseteq
\ell_\infty(\Gamma)$ for convenient index set $\Gamma$, such that
$X/Y$ is non-atomic.
\item All minimal tripotents of $X$ are diagonalizing.
\end{equi}
\end{theorem}

\begin{proof}
$(i)\Rightarrow(ii)$. This is the first part of
Lemma~\ref{lemma-LMP}.

$(ii)\Rightarrow (iii)$. Let $K$ be an elementary triple ideal of
$X$. Since $K$ is an $M$-ideal of $X$, it is clear that condition
$(ii)$ goes down to $K$. Since $K$ is an Asplund space,
Corollary~\ref{coro-char-Asplund} gives us that $K$ has the
alternative Daugavet property, so $K=\mathbb{C}$ by
Proposition~\ref{Cartan}.

$(iii)\Rightarrow (iv)$. The $JBW^*$-triple $X^{**}$ decomposed
into its atomic and non-atomic part as
$X^{**}=\mathcal{A}\oplus_\infty \mathcal{N}$. Let us consider
$Y=\mathcal{A}\cap X$. By \cite[Proposition~3.7]{BC91MZ}, $Y$ is
an Asplund space and $X/Y$ is non-atomic. Since $K_0(Y)=K_0(X)$
(see \cite[Corollary~3.5]{BC92Pacific}) and $K_0(X)=c_0(\Gamma)$,
Proposition~4.4 of \cite{BC92Pacific} shows that
$$
c_0(\Gamma)\subset Y \subset c_0(\Gamma)^{**}=\ell_\infty(\Gamma).
$$

$(iv)\Rightarrow(i)$. On one hand, $X/Y$ has the alternative
Daugavet property since it has the Daugavet property. On the other
hand, a sight to \cite[Remark~8]{FMP} gives us that $n(Y)=1$ and
thus, $Y$ has the alternative Daugavet property. Since $Y$ is an
$M$-ideal, Proposition~3.4 of \cite{MaOi} gives us that $X$ has
the alternative Daugavet property.

$(iii)\Leftrightarrow(v)$. On one hand, a minimal tripotent of $X$
is diagonalizing if and only if it is diagonalizing in the unique
Cartan factor of $X^{**}$ containing it. On the other hand, a
Cartan factor has a diagonalizing minimal tripotent (if and) only
if it equal to $\mathbb{C}$; indeed, if a Cartan factor $C$ has a
minimal diagonalizing tripotent $u$, then $C=C_0(u)\oplus
C_2(u)=C_0(u)\oplus \mathbb{C} u$; but the above direct sum is an
$\ell_\infty$-sum (see \cite[Lemma~1.3]{FR85}) so, by definition
of factor, $C_0(u)=0$.
\end{proof}

The equivalence $(i)\Leftrightarrow (v)$ of the above theorem
appear in \cite[Theorem~4.5]{MaOi} with a different proof..

\begin{remark}
It is worth mentioning that condition $(ii)$ of
Theorem~\ref{th:comm+DPforC} does not necessarily imply the
alternative Daugavet property. For instance, let $H$ be the
two-dimensional Hilbert space and let us consider $X=\ell_1(H)$.
Since $X$ has the Radon-Nikod\'{y}m property and $n(X)=n(\ell_1(H))<1$
by \cite[Proposition~1]{MarPay}, $X$ does not have the alternative
Daugavet property (see \cite[Remark~6]{MarPay}). But since the norm
of $X$ is not Fr\'{e}chet-smooth at any point, the unit ball of $X^*$
does not have any $w^*$-strongly exposed point, thus condition
$(ii)$ of Theorem~\ref{th:comm+DPforC} is satisfied.
\end{remark}

\section{$C^*$-algebras and von Neumann preduals}
The last section of the paper is devoted to particularize the
results of the above section to $C^*$-algebras and von Neumann
preduals. Let us introduce the basic definitions and some
preliminary results.

If $X$ is a $C^*$-algebra, we write $Z(X)$ for the center of the
algebra, i.e.\ the subalgebra consisting in those element of $X$
which commutes with all elements of $X$. A \emph{projection} in
$X$ is an element $p\in X$ such that $p^*=p$ and $p^2=p$. The
tripotents of $X$ are the \emph{partial isometries}, i.e.,
elements $u\in X$ satisfying that $uu^*u=u$. It is clear that
projections are partial isometries (and so tripotents), but there
are partial isometries which are not projections. A projection $p$
in $X$ is said to be \emph{atomic} if $p\neq 0$ and $p X
p=\mathbb{C} p$. The $C^*$-algebra $X$ is said to be
\emph{non-atomic} if it does not have any atomic projection. It is
easy to show (see \cite[pp.~170]{MaOi} or the paragraph before
\cite[Corollary~4.4]{BM}, for example) that a partial isometry $u$
in $X$ is minimal if and only if $d=u^*u$ and $r=u u^*$ (the
domain and range projections associated to $u$) are atomic.
Therefore, a $C^*$-algebra is non-atomic (as an algebra) if and
only if it does not have any minimal tripotent, i.e., it is
non-atomic as a $JB^*$-triple.

The following two corollaries particularize the main results of
the paper to the case of $C^*$-algebras. The first one is a direct
consequence of Theorem~\ref{theorem-JBW}; the second one follows
from Theorem~\ref{th:comm+DPforC} with just a little bit of work.

\begin{corollary}
Let $X$ be a von Neumann algebra. Then, the following are
equivalent:
\begin{equi}
\item $X$ has the alternative Daugavet property.
\item $X_*$ has the alternative Daugavet property.
\item $|x(x_*)|=1$ for every $x\in \ext[X]$ and every $x_*\in
\ext[X_*]$.
\item $X=\mathcal{C}\oplus_\infty \mathcal{N}$ where $\mathcal{C}$
is commutative and $\mathcal{N}$ is non-atomic.
\end{equi}
\end{corollary}

\begin{corollary}
Let $X$ be a $C^*$-algebra. Then, the following are equivalent:
\begin{equi}
\item $X$ has the alternative Daugavet property.
\item $|x^{**}(x^*)|=1$ for every $x^{**}\in \ext[X^{**}]$ and every
$w^*$-strongly exposed point $x^*$ of $B_{X^*}$.
\item There exists a two-side commutative ideal $Y$ such that
$X/Y$ is non-atomic.
\item $K_0(X)$ is isometric to $c_0(\Gamma)$.
\item $K_0(X)$ is commutative.
\item All atomic projections in $X$ are central.
\item $K_0(X)\subseteq Z(X)$.
\end{equi}
\end{corollary}

\begin{proof}
$(i)\Leftrightarrow (ii)\Leftrightarrow (iii)\Leftrightarrow(iv)$
follows from Theorem~\ref{th:comm+DPforC}.

$(iv)\Leftrightarrow (v)$. We have that $K_0(X) =
\left[\oplus_{\gamma\in\Gamma} K(H_\gamma)\right]_{c_0}$ for
convenient index set $\Gamma$ and convenient family of Hilbert
spaces $\{H_\gamma\}$. Thus, $K_0(X)$ is commutative if and only
if $\dim(H_\gamma)\leqslant 1$ for every $\gamma\in \Gamma$.

$(v)\Rightarrow (vi)$. First, we observe that every atomic
projection of $X$ is contained in $K_0(X)$. Now, we take an atomic
projection $p\in X$ and we claim that $p\in Z(X)$. Indeed, fix
$x\in X$ and observe that $p, px, xp\in K_0(X)$ (since it is an
ideal of $X$), so
$$px=p(px)=(px)p=p(xp)=(xp)p=xp.$$

$(vi)\Rightarrow (vii)$. We are going to show that every minimal
partial isometry on $X$ is central, which implies $K_0(X)\subseteq
Z(X)$. Indeed, let $u\in X$ be a minimal partial isometry. We
write $d=u^*u$ and $r=u u^*$ for the domain and range projections
of $u$, which are atomic and so, $d,r\in Z(X)$. Then
$$r=r^2=(uu^*)(uu^*)=
u\big((u^*u)u^*\big)=u\big(u^*(u^*u)\big)=u(u^*u^*)u=\lambda u,$$
where the last equality comes from the fact that $u$ is minimal.
Since $r\neq 0$ (otherwise $u=ru=0$), $u\in \text{span}
(r)\subseteq Z(X)$.

$(vii)\Rightarrow (v)$ is immediate.
\end{proof}

We finish the paper by compiling some known results concerning
$C^*$-algebras and von Neumann preduals having numerical index
$1$. Since we have not found explicitly some of these results in
the literature, and we can now relate them to the alternative
Daugavet property, we include them here for the sake of
completeness.

$C^*$-algebras with numerical index $1$ have been characterized in
\cite{Hur} (see also \cite[pp.~202]{KMR}) as the commutative ones.
With not much work, we can obtain the following result.

\begin{prop}\label{prop-C*-nx1}
Let $X$ be a $C^*$-algebra. Then, the following are equivalent:
\begin{equi}
\item $n(X^*)=1$.
\item $X^*$ has the alternative Daugavet property.
\item $|x^{**}(x^*)|=1$ for every $x^{**}\in \ext[X^{**}]$ and every
$x^*\in \ext[X^*]$.
\item $n(X)=1$.
\item $X$ is commutative.
\end{equi}
\end{prop}

\begin{proof}
$(i)\Rightarrow (ii)$ is clear. $(ii)\Rightarrow (iii)$ is
Corollary~\ref{coro-dualADPimplicanx1}. $(iii)\Rightarrow (iv)$ is
Lemma~\ref{lemma-sufficient-for-nx=1}. $(iv)\Rightarrow (v)$ is
the already mentioned result of \cite{Hur}. Finally,
$(v)\Rightarrow (i)$ follows from the fact that $X^{**}$ is also
commutative and so, $n(X^*)\geqslant n(X^{**})=1$.
\end{proof}

\begin{remark}
In \cite[Example~4.4]{MaOi} it is shown an example of a Banach
spaces with the alternative Daugavet property whose dual does not
share this property. The above proposition gives a lot of examples
of this kind: {\slshape every non-atomic and non-commutative
$C^*$-algebra.} Actually, the example given there fits this
scheme.
\end{remark}

For von Neumann preduals, we have the following.

\begin{prop}\label{prop-vN-nx1}
Let $X$ be a von Neumann algebra. Then, the following are
equivalent:
\begin{equi}
\item $n(X)=1$.
\item $|x^*(x)|=1$ for every $x^*\in \ext[X^*]$ and every
$x\in\ext[X]$.
\item $n(X_*)=1$.
\end{equi}
\end{prop}

\begin{proof}
$(i)\Rightarrow (ii)$ is immediate, since $X$ is commutative and
hence, a $C(K)$ space.

$(ii)\Rightarrow (iii)$ follows from
Lemma~\ref{lemma-sufficient-for-nx=1} applied to $X_*$.

$(iii)\Rightarrow (i)$ is done in \cite[Proposition~1.4]{KMR}.
\end{proof}

Let us finish the paper by raising the following question:
{\slshape is it true that, for a Banach space $Y$, $n(Y^*)=1$
whenever $n(Y)=1$?} We have seen (Propositions \ref{prop-C*-nx1}
and \ref{prop-vN-nx1} above) that this is the case when either $Y$
is a $C^*$-algebra or $Y$ is the predual of a von Neumann algebra.
We do not know the answer even in the case when $Y$ is a
$JB^*$-triple or when $Y$ is the predual of a $JBW^*$-triple.

\vspace{1cm}

\textbf{Acknowledgment:} The author would like to express his
gratitude to \'{A}ngel Rodr\'{\i}guez Palacios for his valuable suggestions
and to Armando R.\ Villena for fruitful conversations about the
subject of this paper.

\end{document}